\documentclass[12pt,oneside,reqno]{article}
\usepackage{amsfonts}
\usepackage{cite}
\usepackage[mathscr]{eucal}
\usepackage{graphicx}
\usepackage{comment}
\usepackage{subcaption}
\newtheorem{theorem}{Theorem}[section]

\newcommand*\C{\mathbb{C}}
\newcommand*\conj[1]{\overline{#1}}
\begin{document}
\title{Cubic Self-inversive Polynomials whose roots envelope conics}
\author{William Calbeck}

\begin{abstract}We study a one parameter family of cubic self-inversive polynomials that "envelope" conic sections in the following sense.  Provided the three roots of the polynomial lie on the unit circle, when you draw the triangle connecting the roots, the sides of the triangle, or extensions thereof, will all be tangent to the same ellipse or hyperbola independent of the parameter.   
\end{abstract}
\maketitle

\section{Introduction}
Let $P(z)$ be a complex polynomial of degree $n$. We say that $P(z)$ is {\it self-inversive } if there exists a complex number $\lambda$ with $|\lambda|=1$ such that 
\begin{equation}\label{sidef}
P(z)=\lambda z^n \conj{P}\left(\frac{1}{z}\right)
\end{equation}
If $P(z)$ has real coefficients and $\lambda=1$ then $P(z)$ is also
called self-reciprocal or palindromic, but in general if $P(z)=\sum_{k=0}^n a_kz^k$ then $P(z)$ is self-inversive iff $a_k=\lambda \conj{a}_{n-k}$ for $i=0,\ldots,n$.  If $z$ is a zero
of $P(z)$ then, because of (\ref{sidef}), so is $1/\conj{z}$.  This property makes self-inversive polynomials good candidates for polynomials whose zeros all lie on the unit circle in $\C$ since a
zero is either on the unit circle or occurs as one of a reflective pair of zeros across the unit circle.  In fact, if a polynomial has all its zeros on the unit circle, then it must be self-inversive (easy exercise), but in general self-inversive polynomials will not have all their zeros on the unit circle, however if the self-inversive polynomial has odd degree, then at least one zero will be on the unit circle.  
\medskip
The following third degree polynomial is self-inversive
\begin{equation}\label{pe}
P_\lambda(z)=z[(z-a)(z-b)]-\lambda[(1-\conj{a}z)(1-\conj{b}z)]
\end{equation}\label{P} where $a,\ b$, and $\lambda$ are in $\C$ and $|\lambda|=1$.  Polynomials of this form are discussed in \cite{lalin13}.  If $a$ and $b$ are inside the unit circle then $P_\lambda(z)$ will have its three zeros on the unit circle.  This type of polynomial, disguised as a Blaschke product, also appeared in \cite{daepp02}. In that article, instead of solving $P_\lambda(z)=0$ they solved, equivalently, $B(z)=\lambda$, where
\begin{equation}
B(z)=z\left(\frac{z-a}{1-\conj{a}z}\right)\left(\frac{z-b}{1-\conj{b}z}\right)
\end{equation}
After arguing that this equation had three distinct solutions on the 
unit circle for each $\lambda$ the authors proceeded to "run the program", so to speak. For each of a spread of $\lambda$ values on the unit circle, a triangle connecting the solutions of $B(z)=\lambda$ was drawn and an ellipse inscribed inside all of the triangles appeared (see figure \ref{fig:ab-inside}).  The authors then went on to prove that it was, in fact,  an ellipse with foci at $a$ and $b$ and found its equation (see below).  If the reader is familiar with Steiner ellipses then be warned that this ellipse may not be the Steiner ellipse of any of the triangles, but may be the Steiner ellipse of one of them for some values of $a$ and $b$ (see \cite{skubak09}).  

\medskip
The purpose of this article is to go beyond what was done in \cite{daepp02}, by allowing $a$ and$/$or $b$ to lie outside the unit circle.  It turns out that as long as the roots stay on the unit circle, the proof given in \cite{daepp02} goes through with only slight modifications, thus yielding an ellipse or hyperbola that is only partially enveloped by the lines connecting the zeros of $P(z)$.  The part that is enveloped is the part passing through the unit disk.
\medskip

\section{Main Result}

\begin{theorem} Let $a,b\in\C$ and $P_\lambda(z)$ as in (\ref{P}). Assume the three roots, $z_1, z_2, z_3$ of $P_\lambda(z)$ are distinct and lie on the unit circle, then each of the three lines connecting two of the roots  , say $z_1$ and $z_2$, will be tangent to the following conic section 

\begin{equation}\label{ce}
||z-a|\pm |z-b||=|1-\conj{a}b|
\end{equation}
at the point $\zeta_3=( m_1 z_2+m_2 z_3)/( m_1+m_2)$ (with similar formulas for $\zeta_1$ and $\zeta_2$) where $m_1, m_2, m_3$  are real numbers and satisfy
\begin{equation}\label{pf}
\frac{(z-a)(z-b)}{(z-z_1)(z-z_2)(z-z_3)}=\frac{m_1}{z-z_1}+\frac{m_2}{z-z_2}+\frac{m_3}{z-z_3}
\end{equation}
The conic equation (\ref{ce}) will be an ellipse (with the "$+$" sign) if both $a$ and $b$ lie inside the unit circle or they both lie outside the unit circle.  If one of $a$ or $b$ lies inside and the other outside the unit circle then the conic equation (with the "$-$" sign) will be a hyperbola.
\end{theorem}

It is not unusual to have conic sections that are tangent to the three sides of a triangle, so what is significant about this result, is that as the parameter $\lambda$ moves around the unit circle producing different self-inversive polynomials $P_\lambda(z)$ with different roots $z_1, z_2, z_3$ the sides of the varying triangles will all be tangent to the same conic section given by equation (\ref{ce}) which only depends on the initial choice of $a$ and $b$.  In this way, the part of the conic section that lies withing the unit disk, is enveloped by the sides of the triangles.
\medskip
\subsection{Examples}
\begin{figure}
    \begin{center}
    \includegraphics[width=0.4\textwidth]{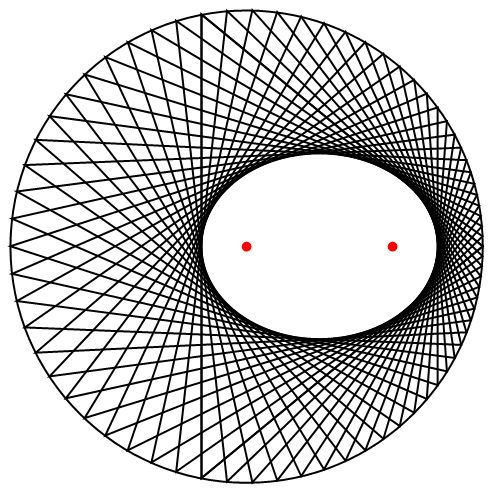}
    \end{center}
    \caption{$P_\lambda(z)=z^3-ez^2+e\lambda z-\lambda$ with $e=0.618$}
    \label{fig:ab-inside}
\end{figure}

Our first example is the case when $a$ and $b$ both lie inside the unit disk as was covered in \cite{daepp02}. In this case, we have an ellipse that lies entirely inside the unit disk.  The ellipse, given by equation (\ref{ce}), has eccentricity $e=|b-a|/|1-\conj{a}b|$. For our example (figure \ref{fig:ab-inside}), to show that any eccentricity is possible, we choose in the simplest case $a=0$ and $b=e$ with $e=0.618$. Recall that, for an ellipse, $e$ ranges between $0$ and $1$.  In this case the enveloping polynomial has the form
\begin{equation}
P_\lambda(z)=z^3-ez^2+e\lambda z-\lambda
\end{equation}

\medskip
\begin{figure}
    \centering
	\begin{subfigure}[b]{0.4\textwidth}
    		\includegraphics[width=\textwidth]{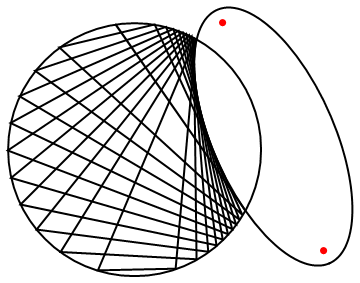}
        	\caption{$a=0.7+i$ and $b=1.5-0.8i$}
    		\label{fig:ab-outside}
	\end{subfigure}
	\begin{subfigure}[b]{0.4\textwidth}
		\includegraphics[width=\textwidth]{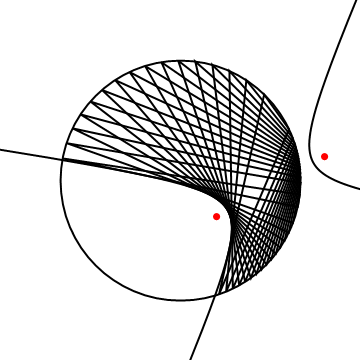}
    		\caption{$a=0.3-0.3I$ and $b=1.2+0.2i$}
    		\label{fig:ab-split}
	\end{subfigure}
	\caption{}
\end{figure}

\medskip
For our second example we choose $a=0.7+i$ and $b=1.5-0.8i$ which both lie outside the unit disk.  Equation (\ref{ce}) still produces an ellipse, but the three roots of $P_\lambda(z)$ will all lie on the unit circle only for certain values of $\lambda$.  We discuss how to find the "good" $\lambda$ values later but it is only for these values that we draw the triangles in figure \ref{fig:ab-outside}.  Note, that one side of each triangle directly contacts the ellipse, but the other sides, if extended, will also make tangential contact with other parts of the ellipse.  This is also true in the last example with $a=0.3-0.3I$ and $b=1.2+0.2i$. In this case, $a$ lies inside the unit disk, $b$ outside, and equation (\ref{ce}) yields a hyperbola, see figure \ref{fig:ab-split}. 
\medskip

\medskip

\subsection{Use of previous results}
Proofs that triangles contain inscribed ellipses, according to Marden \cite{marden45}, date back to Siebeck \cite{siebeck1864}, B{\^o}cher \cite{bocher1892}, Grace \cite{grace01}, Linfield \cite{linfield20}, Marden himself, and more recently Kalman \cite{kalman08} and Badertscher \cite{badertscher14} which deal with the Steiner ellipse, which is the unique ellipse inscribed in a triangle that is tangent at the midpoints of each of the three sides.  As mentioned before, most of the ellipses contructed above are not Steiner ellipses for any of the enveloping triangles. 

\medskip
Let us fix $\lambda$ with $|\lambda|=1$ for the moment and concentrate on just one of the triangles formed by connecting the zeros, which we will call $z_1, z_2$ and $z_3$, of the polynomial $P_\lambda(z)$. \cite{marden45} and \cite{linfield20} tell us that if $m_1,m_2,m_3$ are real and satisfy (\ref{pf})
then there is either an ellipse (if $m_1m_2m_3>0$) or a hyperbola (if $m_1m_2m_3<0$) with foci at $a$ and $b$ that touches each of the sides of the triangle $z_1z_2z_3$.  Where it touches the side $z_iz_k$, it divides it in the ratio $m_i:m_k$.  If one or more of the $m_j$'s are negative then this may occur outside the triangle on the extension of the appropriate side.  The Steiner ellipse occurs when  $a$ and $b$ happen to coincide with the roots of the derivative, $P_\lambda'(z)$.  In this case, $m_1=m_2=m_3=1/3$. This seems to be the result we are looking for, but it doesn't tell us that the ellipse (or hyperbola) is the same ellipse (or hyperbola) for different values of $\lambda$ and it doesn't provide us with the equation of the conic section.  So we take the direct approach as in \cite{daepp02} and show the points of contact, $\zeta_1,\zeta_2,\zeta_3$ satisfy equation (\ref{ce}), independent of $\lambda$.

We follow closely along the lines of the proof in \cite{daepp02} but differ slightly because no mention is made of Blaschke products and we also allow $a$ and/or $b$ to lie outside the unit disk.

\subsection{proof} 
Let $a,b,\lambda\in\C$ with $a,b$ fixed and $|\lambda|=1$.  Let $P_\lambda(z)$ be as in (\ref{pe}) with roots $z_1,z_2,z_3$ which we assume all lie on the unit circle.  Let $m_1,m_2,m_3$ be as in (\ref{pf}).  By comparing coefficients of $z^2$ in (\ref{pf}), it is clear that $m_1+m_2+m_3=1$, and if we multiply
both sides of (\ref{pf}) by $z-z_3$, take the limit as $z\to z_3$, and flip we get
\begin{eqnarray}\label{le}
\frac{1}{m_3} &=& \frac{P'_\lambda(z_3)}{(z_3-a)(z_3-b)} \nonumber\\
&=&\frac{(z_3-a)(z_3-b)+z_3(z_3-b)+z_3(z_3-a)+ \lambda\left[\conj{a}(1-\conj{b}z_3)-\conj{b}(1-\conj{a}z_3)\right]}{(z_3-a)(z_3-b)} \nonumber\\
&=& 1+\frac{z_3}{z_3-a}+\frac{\conj{a}}{\conj{z_3}-\conj{a}}+\frac{z_3}{z_3-b}+ \frac{\conj{b}}{\conj{z_3}-\conj{b}} \nonumber\\
&=& 1+\frac{1-|a|^2}{|z_3-a|^2}+\frac{1-|b|^2}{|z_3-b|^2}
\end{eqnarray}
where we have differentiated (\ref{pe}) and used the fact that $|z_3|=1$ to simplify some of the expressions.  Similar formulas hold for $m_1$ and $m_2$.  These formulas imply that the $m_i$'s are all real, but not necessarily positive.  We again make use of the partial fraction equation (\ref{pf}) by letting $z=a$ to obtain
\begin{eqnarray}\label{za}
0&=&\frac{m_1}{a-z_1}+\frac{m_2}{a-z_2}+\frac{m_3}{a-z_3}\nonumber\\
&=&\frac{m_3}{a-z_3}+\frac{(m_1+m_2)a-(m_1z_2+m_2z_1)}{(a-z_1)(a-z_2)}\nonumber\\
&=&\frac{m_3}{a-z_3}+\frac{(m_1+m_2)(a-\zeta_3)}{(a-z_1)(a-z_2)}\nonumber\\
&=&\frac{m_3}{a-z_3}+\frac{(1-m_3)(a-\zeta_3)}{(a-z_1)(a-z_2)}
\end{eqnarray}
where we have used our definition that $\zeta_3=(m_1z_2+m_2z_1)/(m_1+m_2)$ and $m_1+m_2+m_3=1$.  We can now solve this equation for $\zeta_3-a$ to obtain
\begin{equation}\label{zetaa}
\zeta_3-a=\frac{m_3}{1-m_3}\cdot\frac{(a-z_1)(a-z_2)}{a-z_3}
\end{equation}
in a similar manner using $z=b$ in equation (\ref{pf}) we obtain
\begin{equation}\label{zetab}
\zeta_3-b=\frac{m_3}{1-m_3}\cdot\frac{(b-z_1)(b-z_2)}{b-z_3}
\end{equation}
\medskip
Now we equate the two forms of $P_\lambda(z)$
\begin{displaymath}
z[(z-a)(z-b)]-\lambda[(1-\conj{a}z)(1-\conj{b}z)]=(z-z_1)(z-z_2)(z-z_3)
\end{displaymath}
plug in $z=a$ and manipulate to yield
\begin{displaymath}
\frac{(a-z_1)(a-z_2)}{a-z_3}=\frac{-\lambda(1-|a|^2)(1-\conj{b}a)}{(a-z_3)^2}
\end{displaymath}
which can then be substituted into the right side of (\ref{zetaa}) to give
\begin{equation}\label{zetaa2}
\zeta_3-a=\frac{m_3}{1-m_3}\left[-\frac{\lambda(1-|a|^2)(1-\conj{b}a)}{(a-z_3)^2}\right]
\end{equation}
Similarly
\begin{equation}\label{zetab2}
\zeta_3-b=\frac{m_3}{1-m_3}\left[-\frac{\lambda(1-|b|^2)(1-\conj{a}b)}{(b-z_3)^2}\right]
\end{equation}
Taking absolute values and using the fact that $|\lambda|=1$ and $|1-\conj{a}b|=|1-\conj{b}a|$ for any $a,b\in\C$, we have
\begin{equation}\label{zetaabs}
\pm |\zeta_3-a|\pm|\zeta_3-b|=\left|\frac{m_3}{1-m_3}\right|\cdot |1-\conj{a}b|\cdot\left[\frac{\pm|1-|a|^2|}{|a-z_3|^2}+\frac{\pm|1-|b|^2|}{|b-z_3|^2}\right] 
\end{equation}
The choice of the two $\pm$'s will determine if we have an ellipse or a hyperbola.  We choose the $\pm$'s so that $\pm|1-|a|^2|=1-|a|^2$ and $\pm|1-|b|^2|=1-|b|^2$.  In particular, if both $|a|\leq 1$ and $|b|\leq 1$ then we will choose both $+$'s, if both $|a|>1$ and $|b|>1$ then two $-$'s, otherwise one will be $+$ and one $-$.  In this way we can now use (\ref{le}) to obtain
\begin{equation}\label{finalzeta}
||\zeta_3-a|\pm|\zeta_3-b||=\left|\frac{m_3}{1-m_3}\right|\cdot |1-\conj{a}b|\cdot\left|\frac{1}{m_3}-1\right|=|1-\conj{a}b|
\end{equation}

We have thus shown that $\zeta_3$ satisfies equation (\ref{ce}), so the line passing through $z_1$ and $z_2$ intersects the ellipse or hyperbola at the point $\zeta_3=(m_1z_2+m_2z_1)/(m_1+m_2)$. We make one final use of equation (\ref{pf}) to show that this line is tangent to the ellipse or hyperbola using the well known property of tangent lines to conics, in particular that a line is tangent to an ellipse with foci $a$ and $b$ at a point $\zeta_3$ if it makes equal angles with the lines connecting $\zeta_3$ to $a$ and $\zeta_3$ to $b$.  Plugging $\zeta_3$ into (\ref{pf}) we have

\begin{equation}
\frac{(\zeta_3-a)(\zeta_3-b)}{(\zeta_3-z_1)(\zeta_3-z_2)(\zeta_3-z_3)}=\frac{m_1}{\zeta_3-z_1}+\frac{m_2}{\zeta_3-z_2}+\frac{m_3}{\zeta_3-z_3}
\end{equation}

but
\begin{displaymath}
\frac{m_1}{\zeta_3-z_1}+\frac{m_2}{\zeta_3-z_2}=0
\end{displaymath}

so
\begin{displaymath}
\frac{(\zeta_3-a)(\zeta_3-b)}{(\zeta_3-z_1)(\zeta_3-z_2)}=m_3
\end{displaymath}

Using the complex "$\arg$" function to measure angles and the fact that $m_3$ is real we have
\begin{displaymath}
\arg\left(\frac{\zeta_3-a}{\zeta_3-z_1}\right)+\arg\left(\frac{\zeta_3-b}{\zeta_3-z_2}\right)=0\ {\rm or\ }\pi
\end{displaymath}

If $m_3>0$ then $\zeta_3$ lies between $z_1$ and $z_2$ and the above is "$=0$", otherwise $m_3<0$ and $\zeta_3$ lies on the extension of the line $z_1z_2$ and the above is "$=\pi$". In either case, the line through $z_1$ and $z_2$ intersects the conic section at $\zeta_3$ and makes equal angles there with the lines from $\zeta_3$ to $a$ and $b$ and so it is tangent, or possibly normal, to the conic at $\zeta_3$.  Normalcy can be ruled out by appealing to \cite{linfield20} but we provide the following "sketchy" rationale.  Consider the three examples provided above (or any other examples of the readers choosing).  We have proven that the sides of the triangles intersect the conic sections either tangentially or normally, but the visual evidence allows us to conclude that they intersect tangentially.  We can also have a computer verify that the lines are tangential.  That takes care of at least one example of $a/b$ pairs  for each of the three types.  
Now, because the zeros of a polynomial vary continuously with respect to the coefficients, we can move $a$ and $b$ around and the lines forming the sides of the triangles connecting the zeros  of the polynomials will change directions in a smooth manner.  There will be no sudden $90^\circ$ degree change in direction and so tangential lines will remain tangential. However, this argument breaks down if two of the zeros, say $z_1$ and $z_2$, on the unit circle, were to collide and then abruptly take off in a normal direction to the unit circle.  We mention this because this, in fact, is what happens for fixed $a$ and $b$ as $\lambda$ approaches the limits of the "good" values.  To sucessfully navigate from the examples above to any other similar $a/b$ pair of values, while keeping the $z_i$ distinct, may require adjusting the $\lambda$ value or moving $a$ and/or $b$ closer to the unit circle to expand the good $\lambda$ range.  The isolated points that produce degenerate cases can thus be avoided.  

Computations similar to the one above show that $\zeta_1$ (on the line $z_2z_3$) and $\zeta_2$ (on the line $z_1z_3$) are also points of tangential contact with the conic section.  This completes the proof.

\medskip

\subsection{Finding the good $\lambda$ values}

In order to produce the graphs above we need to find the good $\lambda$ values. We say that $\lambda$ is "good" if the roots $z_1,z_2,z_3$ of $P_\lambda(z)$ are distinct and all lie on the unit circle.  When both $a$ and $b$ lie inside the unit disk, all $|\lambda|=1$ values are good and the ellipse is completely enveloped as in figure \ref{fig:ab-inside}.  So assume $a,b\in\C$ are not both contained in the unit disk. Equation $\ref{ce}$ then produces either an ellipse or hyperbola that merely passes through the unit disk but
is not completely contained therein.  The tangential triangles at the intersection of the conic section with the unit circle are degenerate "triangles" where two of the roots of $P_\lambda(z)$ have converged.  Thus to find the range of good $\lambda$ values we find the $\lambda$ values that produce double roots.  We have done this
computationally solving $P_\lambda(z)=0$ and $P'_\lambda(z)=0$ for $|\lambda|=1$ and assuming all roots are on the unit circle.

\medskip

\section{Conclusion}

We only gave three examples, but any ellipse or hyperbola that passes through the unit disk can be, at least partially, enveloped by a family of self inversive cubic polynomials.  The part that is enveloped is the part that is inside the unit disk.  Other parts of the conic section are also tangent to extensions of sides of the triangle outside the unit circle, but not all of the conic section is enveloped in this way.

\bibliography{selfrecursive}{}
\bibliographystyle{plain}

\bigskip
  \footnotesize

  William Calbeck, \textsc{Department of Mathematics and Physical Sciences, LSU-Alexandria,
    Alexandria, Louisiana 71302}\par\nopagebreak
  \textit{E-mail address}, William Calbeck: \texttt{billc@lsua.edu}
\end{document}